\documentclass{article} 

\usepackage{latexsym}
\usepackage{amsfonts}
\usepackage{graphicx}
\usepackage{subeqn}
\usepackage{cite}
\usepackage{multicol}

\def\BibTeX{{\rm B\kern-.05em{\sc i\kern-.025em b}\kern-.08em
    T\kern-.1667em\lower.7ex\hbox{E}\kern-.125emX}}

\begin{document}

\title{Fast Solution Methods for Fractional Differential Equations in the Modeling of Viscoelastic Materials$^\dag$} 

\author{Kai Diethelm\footnote{FANG, University of Applied Sciences W\"urzburg-Schweinfurt,
Ignaz-Sch\"on-Str.\ 11, 97421 Schweinfurt, Germany, 
Email: kai.diethelm@fhws.de}}

\date{November 22, 2021\footnotetext[2]{%
	IEEE Copyright Notice: \copyright 2021 IEEE. Personal use of this material is permitted. Permission from 
	IEEE must be obtained for all other uses, in any current or future media, including reprinting/republishing 
	this material for advertising or promotional purposes, creating new collective works, for resale or redistribution 
	to servers or lists, or reuse of any copyrighted component of this work in other works.\\
	Accepted for publication in \emph{Proceedings of the 9th International Conference on Systems and Control 
	(ICSC'2021), November 24--26, 2021, Caen, France.}}}

\maketitle

\begin{abstract}
Fractional order models have proven to be a very useful tool for the modeling of the 
mechanical behaviour of viscoelastic materials. Traditional numerical solution methods exhibit 
various undesired properties due to the non-locality of the fractional differential operators, 
in particular regarding the high computational complexity and the high memory requirements.
The infinite state representation is an approach on which one can base numerical methods that 
overcome these obstacles. Such algorithms contain a number of parameters that influence the
final result in nontrivial ways. Based on numerical experiments, we initiate a study leading to
good choices of these parameters.
\end{abstract}

\emph{Keywords:} 
Fractional differential equation, fast solution, diffusive representation

\section{Introduction}
\label{sec:intro}

One of the most commonly suggested and used approaches for formulating the stress-strain relations 
for viscoelastic materials is based on the use of Caputo differential operators of fractional
order, cf., e.g., 
Gerasimov \cite{Ge1948}, Caputo \cite{Ca1967} and Dzherbashyan and his collaborators \cite{DN1968}. 
These operators are defined, for the case $0 < \alpha < 1$ that is of interest here, by
\begin{equation}
  \label{eq:def-cap}
  D_{*0}^\alpha y(t) := D_0^\alpha [y - y(0)](t)
\end{equation}
where $D_0^\alpha$ is the Riemann-Liouville differential operator defined by
\begin{equation}
  \label{eq:def-rld}
  D_0^\alpha y(t) := D J_0^{1-\alpha}y(t).
\end{equation}
Here $D$ is the classical first order differential operator and
\begin{equation}
  \label{eq:def-j}
  J_0^\beta y(t) 
    := \frac1{\Gamma(\beta)} \int_0^t (t-s)^{\beta-1} y(s) \, \mathrm ds,
\end{equation}
for $\beta > 0$, is the Riemann-Liouville integral 
operator of order $\beta$. A detailed discussion of the properties of these
operators may be found in \cite[Chapters 2 and 3]{Di2010}.

A typical material model that includes a few other models as special
cases is the fractional Zener model or fractional standard linear solid,
\begin{equation}
  \label{eq:stress-strain-linear}
  a_0 \sigma(t) + a_1 D_{*0}^\alpha \sigma (t) 
     = m \epsilon(t) + b_1 D_{*0}^\alpha \epsilon(t),
\end{equation}
where $\alpha \in (0,1)$ is a material parameter \cite{MS2011,Me1967,Ma2010}.
Depending on the
type of experiment that is to be described (e.g., a creep test or a
relaxation test), stress
$\sigma$ may be given and strain
$\epsilon$ is the unknown function or vice
versa. In either case, the equation is a special case
of the general mathematical formulation
\begin{subequations}
  \label{eq:ivp}
  \begin{equation}
    \label{eq:fde}
    D^\alpha_{*0} y(t) = f(t, y(t))
  \end{equation}
  where the function $y$ (that typically corresponds to one of the
  functions $\sigma$ or $\epsilon$) is sought on some interval $[0,T]$
  whereas the other function from the set $\{ \sigma, \epsilon\}$ is
  used to define the function $f$ on the right-hand side of the
  fractional differential equation (\ref{eq:fde}).
  Under standard conditions on the function $f$, a unique solution $y$
  can be guaranteed \cite[Theorem 6.5]{Di2010} if this
  differential equation is considered in combination with an initial
  condition of the form\index{initial condition}
  \begin{equation}
    \label{eq:ic}
    y(0) = y_0
  \end{equation}
  with an arbitrary $y_0 \in \mathbb R$.
\end{subequations}

An account of the historical development that has lead to such models
is provided by Rossikhin and Shitikova \cite{Ro2010,RS2007}.

For the numerical solution of such intial value problems, 
many methods have been proposed; cf., e.g., \cite{BDST2016,Di2018,Di2019,Ga2018,LZ2015}.
A common feature of a large part of these methods is that, in order to correctly handle
the non-locality of the fractional differential operators, they require a relatively large amount of
time and/or computer memory \cite{DKLMT2021}.


A possible solution for this challenge is to introduce a so-called \emph{diffusive representation}
(or \emph{infinite state representation}) of the fractional differential operator, cf.\ \cite{Mo1998}. 
This is a representation of the form
\begin{subequations}
	\label{eq:rep-nonclass}
	\begin{equation}
		\label{eq:int-nonclass}
		D_{*0}^\alpha y(t) 
			= \int_0^\infty \phi(w, t) \, \mathrm d w
	\end{equation}
where the integrand $\phi$ 
(whose values $\phi(w, t)$ for $w \in (0, \infty)$ are called the
observed system's infinite states at time $t$) solves the inhomogeneous linear first order 
initial value problem
\begin{equation}
	\label{eq:ivp-nonclass}
	\frac{\partial}{\partial t} \phi (w, t) = h_1(w) \phi(w,t) + h_2(w) y'(t), 
	\quad
	\phi(w, 0) = 0,
\end{equation}
with certain functions $h_1, h_2 : (0, \infty) \to \mathbb R$. Many different specific choices
for these functions are admissible and have been suggested in the literature; cf., e.g.,
\cite{Ba2020,BS2010,Ch2005,Di2008,Di2009,HSL2019,KW2021,McL2018,SG2006,SC2006,TR2002,YA2002,ZCSHN2020}.
\end{subequations}

A common feature of all these special cases of diffusive representations is that using them allows to reduce
the computational complexity of the fractional differential equation solver from the $O(N^2)$ amount of a 
naive implementation of a classical method or the $O(N \log^2 N)$ or $O(N \log N)$ that can be achieved
by suitable modifications \cite{DF2006,FS2001,Ga2018} to just $O(N)$. Moreover, in this way the $O(N)$
memory requirement of the standard algorithms or the $O(\log N)$ demands of some of the modified versions 
mentioned above is reduced to $O(1)$.

While the memory and run time requirements of such diffusive representation based algorithms
are thus immediately evident, their accuracy and their convergence properties are not so obvious.
Indeed, for many of the methods listed above, only a very incomplete error analysis is available. 
In this paper, we therefore aim to start a systematic analysis of the corresponding properties
of at least one such numerical method. We believe that our findings can also be transferred to
many of the other algorithms.

The core of our approach is the observation that any fractional differential equation solver 
needs to be based on some discretization of the
associated differential operator (or the corresponding integral operator if the initial value problem
is reformulated as a Volterra equation before beginning the numerical treatment). In fact, the
approximation quality that this discretization can offer is decisive for the performance 
of the differential equation solver. Therefore, this paper will concentrate on a number of such
discretizations for one specific diffusive representation of the differential operator. 
In view of the construction described in (\ref{eq:rep-nonclass}), all these discretizations
essentially require two components: A numerical integration method for dealing with the 
integral on the right-hand side of (\ref{eq:int-nonclass}) and a numerical solver for the 
first order initial value problem (\ref{eq:ivp-nonclass}). Each of these components 
introduces a certain error into the final result. We shall now investigate the effects of these 
error components.


\section{The Reformulated Infinite State Scheme}
\label{sec:riss}

The concrete approach that we shall discuss is the
relatively recent ``Reformulated Infinite State Scheme" (RISS) 
of Hinze et al.\ \cite{HSL2019}. It is based on the identity
\begin{eqnarray}
	D^\alpha_{*0} y(t) 
		&=& \sin \frac{\alpha \pi}2 \omega^{\alpha-1} y'(t) 
			- \int_0^\infty \mathcal K(\lambda) \frac{\partial}{\partial t} z(\lambda, t) \, \mathrm d \lambda 
				\label{eq:riss-int} \\
		 & & {}    + \cos \frac{\alpha \pi}2 \omega^{\alpha} y(t) 
			- \omega^2 \int_0^\infty \mathcal K(\lambda) Z(\lambda, t) \, \mathrm d \lambda
			\nonumber
\end{eqnarray}
for $t > 0$
(a slight generalization of (\ref{eq:int-nonclass})) where 
\begin{displaymath}
	\omega = \sqrt{\frac{2-\alpha}\alpha}
	\quad \mbox{ and } \quad
	\mathcal K(\lambda) = \frac {\sin \alpha \pi}\pi \frac{\lambda^\alpha}{\lambda^2 + \omega^2}
\end{displaymath}
and where, for $\lambda > 0$, the functions $z$ and $Z$ satisfy the initial value problems
\begin{equation}
	\label{eq:riss-ivp1}
	\frac{\partial}{\partial t} z(\lambda, t) = y'(t) - \lambda z(\lambda, t),
	\quad
	z(\lambda, 0) = 0,
\end{equation}
and
\begin{equation}
	\label{eq:riss-ivp2}
	\frac{\partial}{\partial t} Z(\lambda, t) = y(t) - \lambda Z(\lambda, t),
	\quad
	Z(\lambda, 0) = 0,
\end{equation}
respectively. For the representation (\ref{eq:riss-int}), we shall look at a number of different
concrete discretizations, starting with the original proposal of \cite{HSL2019}. Their suggestion
is to approximate the integrals over the half-line in (\ref{eq:riss-int}) by a truncated compound 
Gaussian quadrature. This means that one introduces some subdivision points 
$0 = \eta_0 < \eta_1 < \ldots < \eta_K$, writes
\begin{equation}
	\label{eq:riss-subdiv}
	\int_0^\infty \mathcal K(\lambda) \zeta(\lambda, t) \, \mathrm d \lambda 
	\approx \sum_{k=1}^K \int_{\eta_{k-1}}^{\eta_k} \mathcal K(\lambda)  \zeta(\lambda, t) \, \mathrm d \lambda 
\end{equation}
(where $\zeta$ may be either $Z$ or $\frac{\partial z}{\partial t}$) and approximately computes each of the 
integrals on the right-hand side of (\ref{eq:riss-subdiv})  
by a $J$-point Gauss-Legendre quadrature. The specific suggestion from
\cite{HSL2019} that we shall follow here first
is to choose $J = 10$, $K = 25$ and $\eta_k = 10^{-5 + 10 (k-1) / (K-1)}$ for $k = 1, 2, \ldots, K$,
so that the tail of
the integration range which is discarded is the interval $[10^{10}, \infty)$.
For the approximate solution of the initial value problems (\ref{eq:riss-ivp1}) and (\ref{eq:riss-ivp2}),
the authors of \cite{HSL2019} suggest to use MATLAB's {\tt ode15s} function. This is a proposal that we
shall not adopt here because {\tt ode15s} may switch between different concrete ODE solvers, making 
a precise analysis virtually impossible. Rather, we follow the advice of \cite{Ba2020,Di2009} in a similar situation
and use an A-stable scheme (which may not necessarily be the case in {\tt ode15s}). More specifically, 
we will test both the backward Euler formula (a first-order method), viz.\
\begin{subequations}
\begin{equation}
	\label{eq:ode-z-bweuler}
	z(\lambda, t_{n}) \approx \frac 1 {1 + h \lambda} \left( z(\lambda, t_{n-1}) + h y'(t_n) \right)
\end{equation}
and
\begin{equation}
	Z(\lambda, t_{n}) \approx \frac 1 {1 + h \lambda} \left( Z(\lambda, t_{n-1}) + h y(t_n) \right)
\end{equation}
\end{subequations}
where $h$ is the step size,
and the trapezoidal method (which is of second order),
\begin{subequations}
\begin{equation}
	\label{eq:ode-z-trap}
	z(\lambda, t_{n}) 
	\approx \frac { z(\lambda, t_{n-1})  \left[ 1 - \frac h 2 \lambda \right] 
						+ \frac h 2 [ y'(t_n) + y'(t_{n-1}) ] } {1 + h \lambda/2}
\end{equation}
and
\begin{equation}
	Z(\lambda, t_{n}) 
	\approx \frac { Z(\lambda, t_{n-1})  \left[ 1 - \frac h 2 \lambda \right] 
						+ \frac h 2 [ y(t_n) + y(t_{n-1}) ] } {1 + h \lambda/2} .
\end{equation}
\end{subequations}

It may be observed that the kernel $\mathcal K$ in the integrals in (\ref{eq:riss-int}) behaves as
$\mathcal K(\lambda) \sim \lambda^\alpha$ for $\lambda \to 0$ and so, since $0 < \alpha < 1$,
it is not differentiable at the origin. This is a situation that the compound Gauss method in its
original form cannot handle
very well. Therefore, we also investigate the behaviour of a slightly modified scheme where we
replace the application of the Gauss-Legendre formulas to the integrals on the right-hand side of 
(\ref{eq:riss-subdiv}) by (a) the standard Clenshaw-Curtis (CC) formula \cite{BP2011} for those subintervals
$[\eta_{k-1}, \eta_k]$ where $\eta_{k-1} \ge 1$ (i.e.\ the subintervals that are sufficiently far away 
from the singular point of the kernel $\mathcal K$), and (b) a weighted Clenshaw-Curtis method with weight
function $\lambda^\alpha$, thus accurately capturing the asymptotics of $\mathcal K(\lambda)$,
for the remaining subintervals.

\section{A Comparison of the Fundamental Variants of the Algorithm}

We start by looking at the run times of the algorithm.
Table \ref{tab:timing1} shows the run times required for a typical example (the results for other examples 
are comparable). In particular, we can see that, as expected  from the construction of the algorithm,
the run times are approximately proportional to $J K h^{-1}$. We can also see that the trapezoidal ODE solver,
having a more complex structure than the backward Euler method (in particular, requiring more
evaluations of the function $y$), usually requires a significantly longer run time (more precisely:
up to almost three times as long) than backward Euler. However, we will see below that this
increased run time often comes with a significant accuracy benefit that more than makes up 
for the higher cost.

\begin{table}[h]
	\caption{\label{tab:timing1}Run times (in seconds) for computation of $D^\alpha_{*0} y(t)$ for $t \in [0,3]$ 
		with $\alpha = 0.4$ and $y(t) = t^{1.6}$ (RISS, exact evaluation of $y'$)}
	\centering
	\addtolength{\tabcolsep}{-0.5mm}
	\begin{tabular}{l|l|r|r|r|r}
		\multicolumn{2}{c|}{\ } & \multicolumn{2}{c|}{backward Euler} & \multicolumn{2}{c}{trapezoidal} \\
		\multicolumn{2}{c|}{\ } & \multicolumn{2}{c|}{ODE solver} & \multicolumn{2}{c}{ODE solver} \\
		\hline
		\multicolumn{2}{c|}{ODE solver time step size $h$} & $10^{-2^{\vphantom{2}}}$ & 
				$10^{-4}$  & $10^{-2}$ & $10^{-4}$ \\
		\hline
		compound & $J = 25$, $K = 10$ & 0.547 & 57.20 & 1.484 & 148.4 \\
		\cline{2-6}
		Gauss       & $J = 10$, $K = 25$ & 0.609 & 58.28 & 1.406 & 138.5 \\
		\cline{2-6}
		quadrature & $J = 15$, $K = 7$ & 0.297 & 30.22 & 0.453 & 46.0 \\
		\hline
		compound & $J = 25$, $K = 10$ & 0.641 & 74.09 & 1.406 & 137.9 \\
		\cline{2-6}
		(weighted) CC & $J = 10$, $K = 25$ & 0.578 & 60.00 & 1.359 & 134.2 \\
		\cline{2-6}
		quadrature & $J = 15$, $K = 7$ & 0.266 & 28.05 & 0.453 & 44.1 
	\end{tabular}
\end{table}

\begin{figure}[h]
	\centering
	\includegraphics[width=0.9\columnwidth]{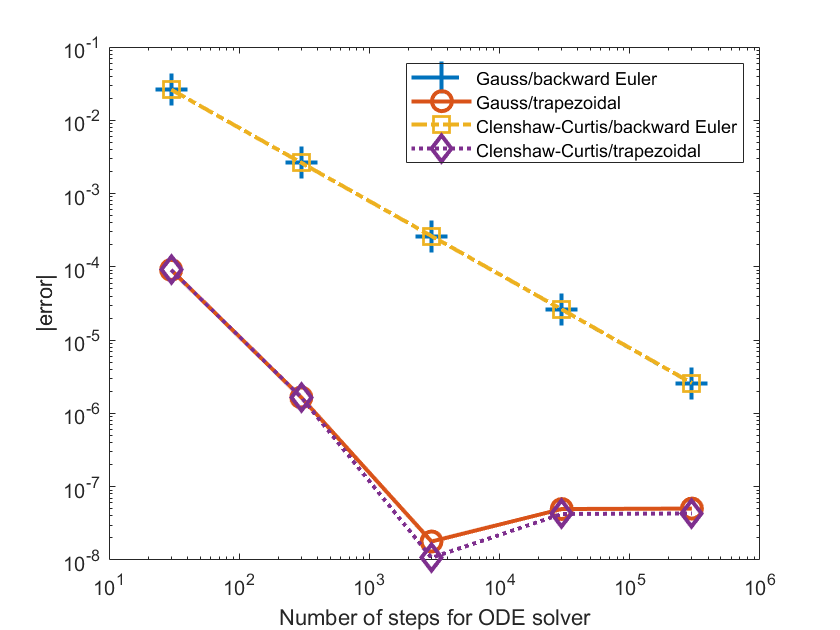}
	\caption{\label{fig:v1-ex1-25-10}Maximal absolute values of approximation errors 
		of RISS methods with different discretizations for computation of $D^{0.4}_{*0} y(t)$
		with $y(t) = t^{1.6}$. In all cases, $J = 25$ and $K = 10$.}
\end{figure}

\begin{figure}[h]
	\centering
	\includegraphics[width=0.9\columnwidth]{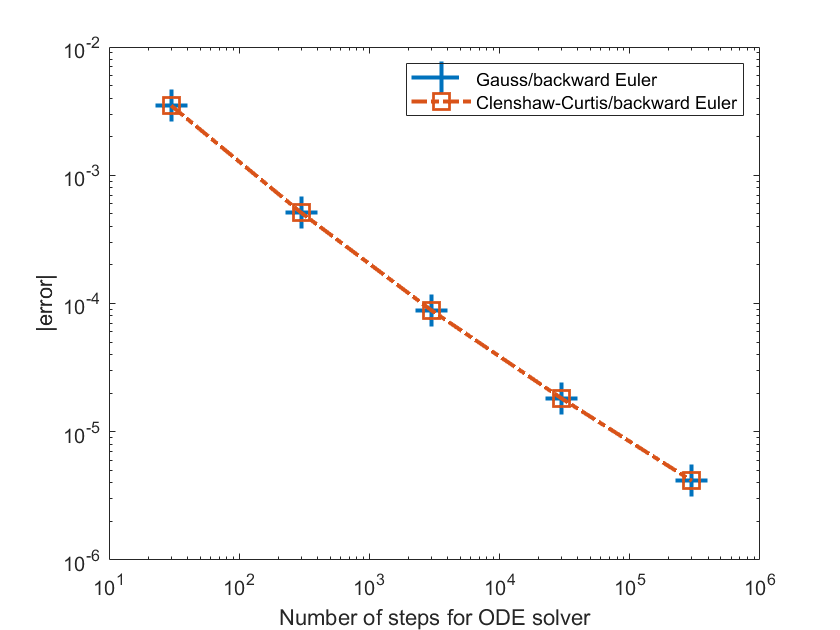}
	\caption{\label{fig:v1-ex3-25-10}Maximal absolute values of approximation errors 
		of RISS methods with different discretizations for computation of $D^{0.4}_{*0} y(t)$
		with $y(t) = t^{0.6}$. In all cases, $J = 25$ and $K = 10$.}
\end{figure}

Regarding the accuracy of the numerical results, we first refer to Figs.\ \ref{fig:v1-ex1-25-10}
and \ref{fig:v1-ex3-25-10}. 
In Fig.\  \ref{fig:v1-ex1-25-10}, we have looked at a rather well behaved example, 
viz.\ the computation of $D^{0.4}_{*0} y(t)$
with the differentiable function $y(t) = t^{1.6}$, for $t \in [0, 3]$. We have used the RISS scheme
with parameters $J = 25$ and $K = 10$ for the numerical quadrature. All possible combinations
of quadrature formula (compound Gauss or compound CC) and ODE solver (backward Euler or
trapezoidal) have been tried with different numbers of time steps in the ODE solver.
Each of the four curves in that figure shows the maximal error over the entire interval
for one of these four combinations. 
Figure \ref{fig:v1-ex3-25-10} shows the analog results for $y(t) = t^{0.6}$ (a much less smooth function,
but nevertheless a function with a behaviour that is typical for solutions to fractional 
differential equations),
with all other parameters being the same as in Fig.\ \ref{fig:v1-ex1-25-10}.

From these two examples, we can already observe a few properties that are also valid in
general:
\begin{itemize}
\item The choice of the quadrature formula (compound Gauss or compound (weighted) 
	Clenshaw-Curtis) only has a negligible effect on the final result. This may be considered
	to be slightly unexpected because, as indicated at the end of Section \ref{sec:riss}, the idea
	behind considering the (weighted) Clenshaw-Curtis scheme was that the Gauss method
	might not perform well. However, this observation has a natural explanation: We do not 
	use the Gauss method in its standard form but a compound Gauss method with a highly graded
	mesh (the mesh grading being described by the choice of the parameters $\eta_k$). Such
	graded meshes indeed are very suitable for treating singularities of the type encountered here,
	and it is the influence of this mesh grading (which is the same for both choices of
	the quadrature formulas) that leads to this small difference between Gauss and Clenshaw-Curtis.
\item The overall error that is displayed in the figure consists of two components due to the
	quadrature formula and the ODE solver, respectively.
\item For the trapezoidal method, a certain saturation is reached when the error reaches
	a level of about $5 \cdot 10^{-8}$. This indicates that the error of the quadrature formula
	(which, as indicated above, is almost exactly the same in all cases) approximately has this value.
	For sufficiently small step sizes, the ODE solver's error becomes so small that the quadrature error
	completely dominates the overall error.
\item In both examples, the overall error decreases to this level as $O(h)$ for the backward Euler method. 
	This matches the well known fact that the backward Euler method is of first order.
\item For the trapezoidal method, the decrease is $O(h^{1.6})$ in the example shown in
	Fig.~\ref{fig:v1-ex1-25-10}. The full second order convergence that
	this method can have under optimal circumstances is not achieved because the function 
	$y \notin C^2$.
\item For the example presented in Fig.\ \ref{fig:v1-ex3-25-10}, we cannot directly apply the trapezoidal method 
	because, as can be seen from (\ref{eq:ode-z-trap}) for $n=1$, the approximate computation of the infinite
	states $z(\lambda, t_1)$ requires the evaluation of $y'(t_0) = y'(0)$ which does not exist. We shall discuss 
	some ways to successfully deal with this challenge in Section \ref{sec:refine}.
\end{itemize}

The plots of Figs.\ \ref{fig:v1-ex1-25-10} and \ref{fig:v1-ex3-25-10} demonstrate the effects of
the choice of the quadrature formula, the ODE solver and the solver's step size as the integration range 
decomposition and the number of quadrature points remain unchanged. The next step of our investigation
is to see what will happen when these latter parameters do change. Some results in this connection
are shown in Table \ref{tab:jk-change}.

\begin{table}[h]
	\caption{\label{tab:jk-change}Absolute values of errors for the example from Fig.\ \ref{fig:v1-ex1-25-10}
		for different choices of $J$ and $K$}
	\centering
	\addtolength{\tabcolsep}{-0.5mm}
	\begin{tabular}{l|c|c|r|r|r|r}
		\multicolumn{3}{c|}{\ } & \multicolumn{2}{c|}{backward Euler} & \multicolumn{2}{c}{trapezoidal} \\
		\multicolumn{3}{c|}{\ } & \multicolumn{2}{c|}{ODE solver} & \multicolumn{2}{c}{ODE solver} \\
		\hline
		\multicolumn{3}{c|}{time step size $h$} & $10^{-2^{\vphantom{2}}}$ & 
				$10^{-4}$  & $10^{-2}$ & $10^{-4}$ \\
		\hline
		& $J$ & $K$ & & & & \\
		\hline
		compound & 25 & 10 & 2.63E$-$3 & 2.63E$-$5 & 1.65E$-$6 & 4.96E$-$8 \\
		\cline{2-7}
		Gauss       & 10 & 25 & 2.63E$-$3  & 2.63E$-$5 & 1.66E$-$6  & 4.43E$-$8 \\
		\cline{2-7}
		quadrature & 15 & 7 & 2.63E$-$3  & 2.15E$-$5 & 3.17E$-$6 & 4.87E$-$6 \\
		\hline
		compound & 25 & 10 & 2.63E$-$3 & 2.63E$-$5 & 1.66E$-$6 & 4.24E$-$8 \\
		\cline{2-7}
		(weighted) & 10 & 25 & 2.63E$-$3 & 2.63E$-$5 & 1.64E$-$6 & 6.35E$-$8 \\
		\cline{2-7}
		CC quadr. & 15 & 7 & 2.83E$-$3 & 2.26E$-$4 & 2.01E$-$4 & 1.99E$-$4 
	\end{tabular}
\end{table}

The results of Table \ref{tab:jk-change} indicate the following facts:
\begin{itemize}
\item As in the figures above, we can see that the backward Euler method induces a relatively large 
	error that dominates the entire process. In other words, if one wants to apply this ODE solver with
	the step sizes used here, it is possible to significantly reduce the values of the parameters $J$
	and $K$, and hence also the computational complexity, without any noticeable adverse effects on the final 
	result.
\item For the case $J = 15$ and $K = 7$, the compound Gauss quadrature produces significantly better results 
	than the compound (weighted) Clenshaw-Curtis method for most choices of the
	ODE solver. In other examples, e.g.\ when the order of the derivative is changed from $0.4$ to 
	$0.82$, the opposite can be observed. For this behaviour, no explanation is currently available.
\end{itemize}

\section{A Refinement}
\label{sec:refine}

As indicated above, it follows from eq.~(\ref{eq:ode-z-trap}) that 
the trapezoidal method can only be used as the ODE solver in the scheme
under consideration if the function $y$ is differentiable throughout the
entire interval of interest and, in particular, at least has one-sided derivatives
at the end points of this interval. However, it is well known, cf.\ \cite[\S6.4]{Di2010} or
\cite{Di2007}, that this assumption is satisfied at the initial point, i.e.\ at the left end point 
of the interval in question, only in exceptional cases. Therefore, if one would like to
work with the trapezoidal method in our context, then one needs to find a way to avoid the
explicit use of the function $y'$ in eq.~(\ref{eq:ode-z-trap}).

The most natural approach to deal with this challenge is to replace the terms $y'(t_n)$ and
$y'(t_{n-1})$ by appropriately chosen finite differences. To this end, we have tested a
procedure that appears to be a rather natural concept. The construction of this
scheme is based on exploiting the following observations:
\begin{itemize}
\item The expressions  $y'(t_n)$ and  $y'(t_{n-1})$
	do not appear in eq.~(\ref{eq:ode-z-trap}) in an independent manner; rather, only their 
	sum is relevant.
\item The reason for using the trapezoidal formula is its higher accuracy when compared to the
	backward Euler method. This will be lost if one uses a very rough approximation of
	the derivative; hence, one should use a finite difference that is second order accurate.
\item In the first step, i.e.\ in the case $n = 1$, one should only use a forward difference
	to approximate the term $y'(t_{n-1}) = y'(0)$
	as otherwise one would need to evaluate $y(-h)$. This clearly should be avoided because
	it cannot be guaranteed that $y$ is defined for negative arguments.
\item For analog reasons, in the last step one should only use a backward difference for 
	$y'(t_n)$.
\item In all other steps, centered differences may be used. 
\end{itemize}

Therefore, in the first attempt we have chosen to use the discretizations
\begin{equation}
	\label{eq:fin-diff}
	y'(t_n) \approx 
		\frac 1 {2h} \times
		\cases{
			(-y(t_2) + 4 y(t_1) - 3 y(t_0))             &  if  $n = 0$, \cr
			(3 y(t_N) - 4 y(t_{N-1}) + y(t_{N-2}))  &  if  $n = N$, \cr
			(y(t_{n+1}) - y(t_{n-1}))                    &  else, \cr
		}
\end{equation}
in each case of which have an error of $O(h^2)$ if $y \in C^2[t_0, t_N]$.
Even though, as mentioned above, this smoothness assumption is not satisfied in all
common applications, we will nevertheless obtain results that so sufficiently
accurate that their error does not dominate the error induced by the trapezoidal 
formula. It thus follows that our modified algorithm can be constructed by replacing the
term $(h/2) [y'(t_n) + y'(t_{n-1})]$ in the numerator of the right-hand side 
of (\ref{eq:ode-z-trap}) by
\[
		\cases{
			y(t_1) - y(t_0)             &  if  $n =1$, \cr
			y(t_N) - y(t_{N-1})      &  if  $n = N$, \cr
			(y(t_{n+1}) + y(t_n) - y(t_{n-1}) - y(t_{n-2}) ) / 4   &  else. \cr
		}
\]
For the sake of consistency and to obtain even more insight, we also 
modify the backward Euler method in a corresponding way. Here it suffices to
replace the term $h y'(t_n)$ in (\ref{eq:ode-z-bweuler}) by the associated
backward difference, i.e.\ by $(y(t_n) - y(t_{n-1}))$, so that this method also 
does not require to evaluate derivatives of $y$ any more.

\begin{figure}[h]
	\centering
	\includegraphics[width=0.9\columnwidth]{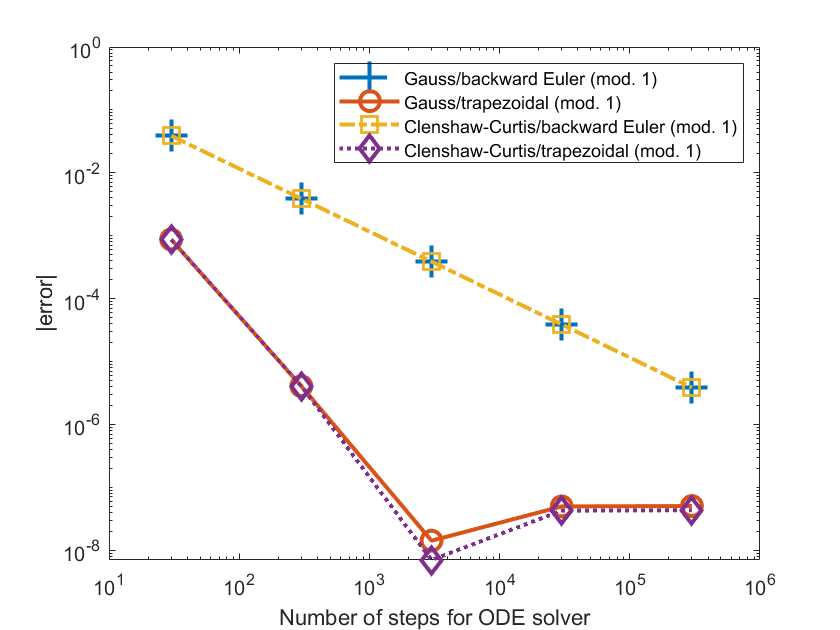}
	\caption{\label{fig:v1-ex1-25-10-riss3}Maximal absolute values of approximation errors 
		of first modified RISS methods with different discretizations 
		for computation of $D^{0.4}_{*0} y(t)$
		with $y(t) = t^{1.6}$. In all cases, $J = 25$ and $K = 10$.}
\end{figure}

Figure \ref{fig:v1-ex1-25-10-riss3} shows the numerical results for the example 
already discussed in Fig.\ \ref{fig:v1-ex1-25-10}, where now the ``standard"
backward Euler and trapezoidal methods given in 
(\ref{eq:ode-z-bweuler}) and (\ref{eq:ode-z-trap}), respectively, 
have been modified in the indicated manner.
It can be seen that the numerical results differ from those obtained with the
unmodified methods only by a very small amount. In particular, all the conclusions
drawn about the behaviour of the methods for this example remain valid.
This indicates that the modifications have introduced an additional error that 
is of the same order as the error that the respective ODE solvers induce anyway,
so that the overall error order remains unchanged.

\begin{figure}[h]
	\centering
	\includegraphics[width=0.9\columnwidth]{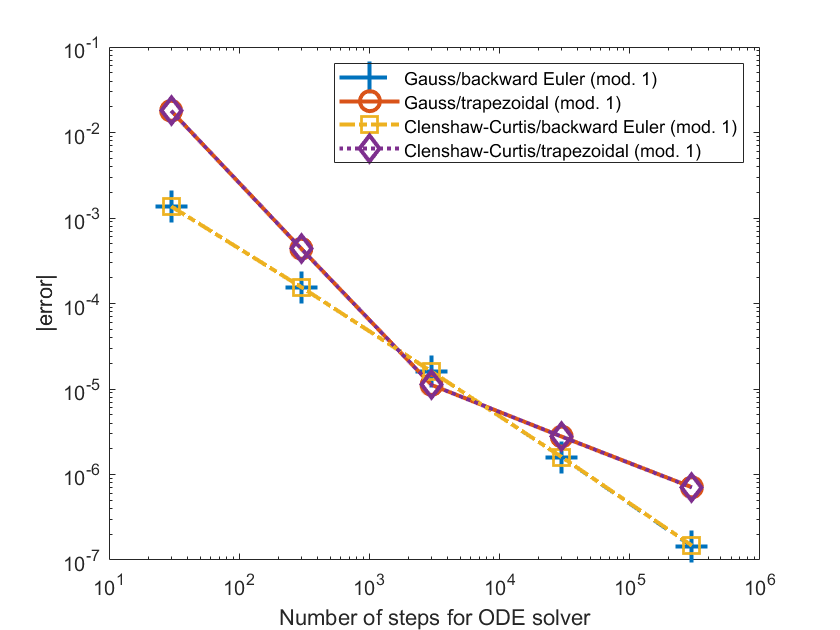}
	\caption{\label{fig:v1-ex3-25-10-riss3}Maximal absolute values of approximation errors 
		of first modified RISS methods with different discretizations 
		for computation of $D^{0.4}_{*0} y(t)$
		with $y(t) = t^{0.6}$. In all cases, $J = 25$ and $K = 10$.}
\end{figure}

Similarly, Fig.\ \ref{fig:v1-ex3-25-10-riss3} shows the numerical results 
of these modified approaches for the example from Fig.\ \ref{fig:v1-ex3-25-10}. 
Here now we can see that our modification admits to use not only the 
backward Euler method (as before) for this example but also the modified trapezoidal 
algorithm which was not possible for the unmodified version. The modified 
backward Euler scheme still exhibits an $O(h)$ convergence towards the saturation
limit induced by the quadrature formula. For the modified trapezoidal formula, the behaviour 
is more complex: The graphs exhibit a clearly recognizable kink. An analysis of the slopes of
the graphs at either side of the kink indicate a convergence behaviour towards the saturation limit
that has the form $c_1 h^{0.6} + c_2 h^{1.6} + o(h^{1.6})$ with some nonzero constants
$c_1$ and $c_2$ having the property that $|c_1|$ is very much smaller than $|c_2|$. 
Therefore, while $h$ is relatively large, $c_1 h^{0.6}$ is still much smaller in absolute
value than $c_2 h^{1.6}$, so the latter contribution dominates the entire error. 
When $h$ becomes smaller, the situation reverses and the graph's slope becomes less steep.

An alternative idea in the context of the trapezoidal method 
is to always use a forward difference,
i.e.\ the method used in the first case of eq.~(\ref{eq:fin-diff}), to approximate the term 
$y'(t_{n-1})$ in the expression $(h/2) [y'(t_n) + y'(t_{n-1})]$ in (\ref{eq:ode-z-trap}) 
and to use the backward difference from the second case in~(\ref{eq:fin-diff}) to replace 
$y'(t_n)$ there. This yields the expression 
\[
	\frac h 2 [y'(t_n) + y'(t_{n-1})] \approx y(t_n) - y(t_{n-1})
\]
that, apart from being of a slightly simpler structure, also has the advantage that
it can be used for all $n$ so that a distinction of cases is no longer necessary. 

Moreover, for the backward Euler method we may also use a second order 
centered difference, i.e.\ we may replace $h y'(t_n)$ by $(y(t_{n+1}) - y(t_{n-1}))/2$.

\begin{figure}[h]
	\centering
	\includegraphics[width=0.9\columnwidth]{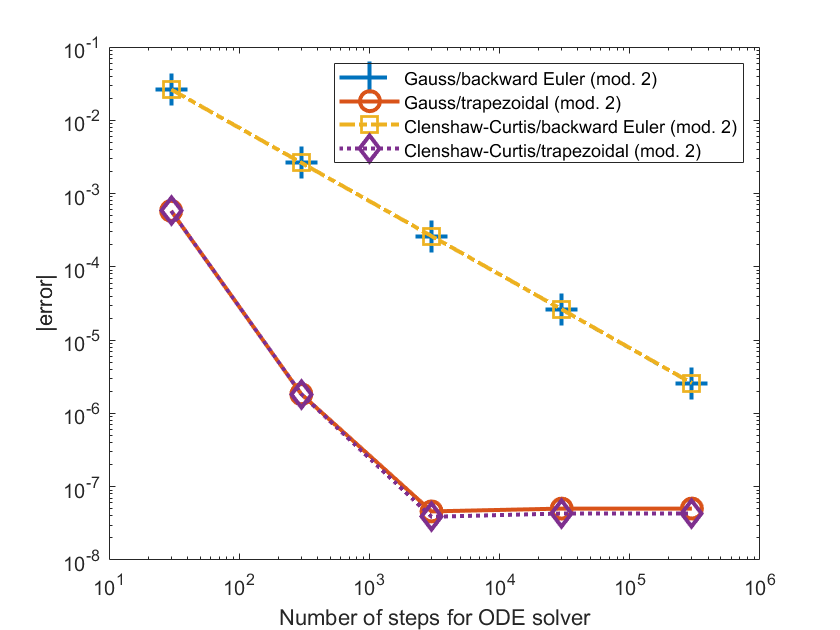}
	\caption{\label{fig:v1-ex1-25-10-riss4}Maximal absolute values of approximation errors 
		of RISS methods using second modification 
		with different discretizations for computation of $D^{0.4}_{*0} y(t)$
		with $y(t) = t^{1.6}$. In all cases, $J = 25$ and $K = 10$.}
\end{figure}

\begin{figure}[h]
	\centering
	\includegraphics[width=0.9\columnwidth]{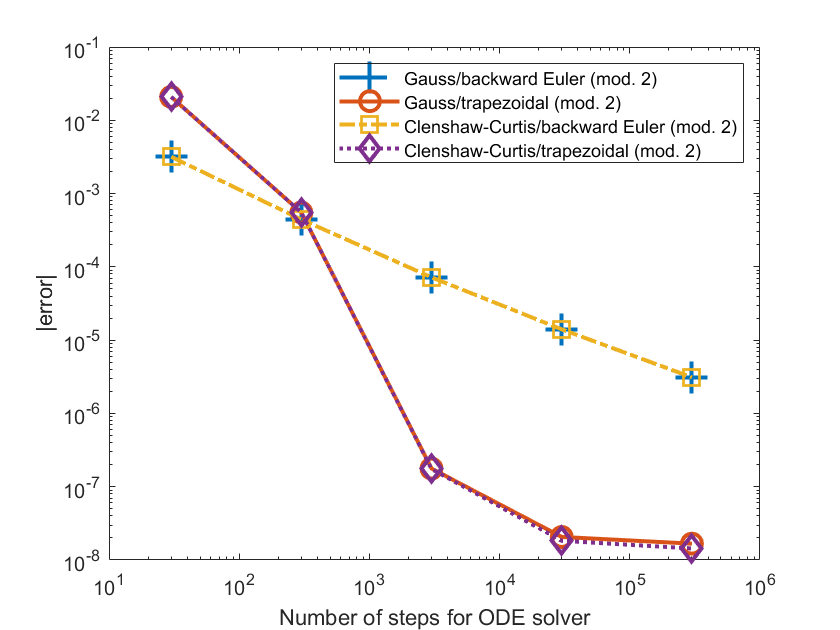}
	\caption{\label{fig:v1-ex3-25-10-riss4}Maximal absolute values of approximation errors 
		of RISS methods using second modification
		with different discretizations for computation of $D^{0.4}_{*0} y(t)$
		with $y(t) = t^{0.6}$. In all cases, $J = 25$ and $K = 10$.}
\end{figure}

The numerical results obtained
for our two examples with these alternative modifications are shown 
in Figs.\ \ref{fig:v1-ex1-25-10-riss4} and \ref{fig:v1-ex3-25-10-riss4}.

For the first example, we once again find almost no difference between this approach
and the original or the first modified version (see Fig.\ \ref{fig:v1-ex1-25-10-riss4}).

The second and less smooth example is shown in Fig.\  \ref{fig:v1-ex3-25-10-riss4}. 
Here, the Euler method with the second modification 
has a convergence order of only $O(h^{0.6})$, 
presumably due to the poor performance of the finite difference approximation to $y'$
for this rather nonsmooth function.
For the second modified trapezoidal method, we seem to have a
simliar behaviour as in the case of the first modified version of the trapezoidal formula,
i.e.\ an error component of the form $c_1 h^{0.6} + c_2 h^{1.6} + o(h^{1.6})$ with 
$0 < |c_1| \ll |c_2|$, although the figures are a bit less clear, so a more detailed 
analysis should be performed.

\section{Conclusion}

Our numerical experiments as well as elementary theoretical considerations
reveal that the error of the approximate calculation of a fractional derivative
using the RISS method comprises three components and thus has the form
\[
	\epsilon_1(\eta_K) + \epsilon_2(J, K, \eta_1, \ldots, \eta_K) + \epsilon_3(h).
\]
Here,
\begin{itemize}
\item $\epsilon_1(\eta_K)$ is the error due to the truncation of integration interval that
	behaves as $\epsilon_3(\eta_K) \to 0$ for $\eta_k \to \infty$,
\item $\epsilon_2(J, K, \eta_1, \ldots, \eta_K)$ is the integration error
	that satisfies\linebreak $\epsilon_2(J, K, \eta_1, \ldots, \eta_K) \to 0$ if $J \cdot K \to \infty$  
	when $\eta_K$ is fixed, and
\item $\epsilon_3(h) = O(h^p)$ is the error due to the ODE solver where
	$p$ depends on choice of the solver and the degree of smoothness of function in question.
\end{itemize}

The indicated modifications of the RISS approach do not lead to any essential changes
in these properties (but they may have an influence on the value of $p$ that determines
the behaviour of $\epsilon_1(h)$).

We believe that it is possible to achieve $\epsilon_3 = 0$ by a proper choice of the quadrature formula,
but the effects of this on $\epsilon_2$ are likely nontrivial.

Analog observations hold for almost all the other numerical methods based on 
such principles mentioned in Section \ref{sec:intro}.

This relatively clear general picture however needs much more detailed additional
investigations to fill in details like a description of the interaction
between the parameters $h$, $J$ and $K$ that might lead to some advice on which 
combinations of values for these parameters are useful in practice.

\bibliographystyle{plain}
\bibliography{icsc}

\end{document}